\documentclass[11pt,fleqn,twoside]{article}
\usepackage{amsfonts,amssymb,latexsym}
\makeatletter
\newcommand{\prava}[1]{\small\it
\begin{flushleft}
Copyright \copyright \ 1999 by  #1
\end{flushleft}}

\newcommand{\name}[1]{\begin{flushleft}
                       \LARGE \bf #1
                       \end{flushleft}\vspace{-3mm}}

\newcommand{\Author}[1]{\begin{flushleft}
                       \it #1 \end{flushleft}}

\newcommand{\Adress}[1]{\begin{flushleft}
                       \it #1 \end{flushleft}}

\newcommand{\Date}[1]{\begin{flushleft}
                      \small  \it #1 \end{flushleft}}

\newcommand{\ehkol}{Author \ name}
\newcommand{\ohkol}{Article \ name}
\renewcommand{\@evenhead}{
\hspace*{-3pt}\raisebox{-15pt}[\headheight][0pt]{\vbox{\hbox to \textwidth 
{\thepage \hfil \ehkol}\vskip4pt \hrule}}}
\renewcommand{\@oddhead}{
\hspace*{-3pt}\raisebox{-15pt}[\headheight][0pt]{\vbox{\hbox to \textwidth 
{\ohkol \hfil \thepage}\vskip4pt\hrule}}}
\renewcommand{\@evenfoot}{}
\renewcommand{\@oddfoot}{}

\newcommand{\be}{\begin{equation}}
\newcommand{\ee}{\end{equation}}
\newcommand{\ba}{\hspace*{-5pt}\begin{array}}
\newcommand{\ea}{\end{array}}

\newcommand{\ds}{\displaystyle}
\makeatother

\newcommand{\zi}{z_i}
\newcommand{\bzi}{\bar z_i}
\newcommand{\bdzi}{\partial_{\bar\zi}}
\newcommand{\dzi}{\partial_{\zi}}
\newcommand{\xii}{\xi_i}
\newcommand{\dxi}{\partial_{\xii}}
\newcommand{\zj}{z_j}
\newcommand{\bzj}{\bar z_j}
\newcommand{\dzj}{\partial_{\zj}}
\newcommand{\bdzj}{\partial_{\bar\zj}}
\newcommand{\xij}{\xi_j}
\newcommand{\dxj}{\partial_{\xij}}
\newcommand{\phii}{\phi_i}
\newcommand{\psii}{\psi_i}
\newcommand{\bph}{\bar\phi}
\newcommand{\bps}{\bar\psi}
\newcommand{\dep}{\delta_{\pi}}

\newcommand{\dpi}{\delta_{\pi}}
\newcommand{\Dpi}{\Delta_{\pi}}
\newcommand{\Lzi}{L_{\zi\dzi}}
\newcommand{\bLzi}{L_{\bzi\bdzi}}
\newcommand{\Lxi}{L_{\xii\dxi}}
\newcommand{\ixi}{i_{\xi\dxi}}
\newcommand{\izi}{i_{\zi\dzi}}
\newcommand{\bizi}{i_{\bzi\bdzi}}
\newcommand{\Lxj}{L_{\xij\dxj}}
\newcommand{\Lzj}{L_{\zj\dzj}}
\newcommand{\bLzj}{L_{\bzj\bdzj}}
\newcommand{\ixj}{i_{\xij\dxj}}
\newcommand{\izj}{i_{\zj\dzj}}
\newcommand{\bizj}{i_{\bzj\bdzj}}
\newcommand{\cn}{A\otimes A^{\otimes n}}

\newtheorem{explanation}{Explanation}
\newtheorem{proposition}{Proposition}

\begin{document}

\thispagestyle{empty}
\setcounter{page}{365}
\renewcommand{\ehkol}{A. Kotov}
\renewcommand{\ohkol}{Poisson Homology of $r$-Matrix Type
   Orbits I: Example of Computation}

\begin{flushleft}
\footnotesize \sf Journal of Nonlinear Mathematical Physics \qquad
1999, V.6, N~4, \pageref{kotov-fp}--\pageref{kotov-lp}.
\hfill {\sc Article}
\end{flushleft}

\vspace{-5mm}

\renewcommand{\footnoterule}{}
{\renewcommand{\thefootnote}{}
 \footnote{\prava{A. Kotov}}}

\name{Poisson Homology of {\mathversion{bold}$r$}-Matrix Type
   Orbits I: Example of Computation}\label{kotov-fp}

\Author{Alexei KOTOV}

\Adress{Mathematical Physics Group, Institute of Theoretical
and Experimental Physics,\\
B. Cheremuchkinskaya ul.25, 117 259 Moscow, Russia\\[1mm]
E-mail:  kotov@pan.teorfys.uu.se, kotov@vitep5.itep.ru}

\Date{Received February 18, 1999; Revised June 14, 1999; Accepted
July 14, 1999}

\begin{abstract}
\noindent
In this paper we consider the Poisson algebraic structure
associated with a classical
$r$-matrix, i.e. with a solution of the modif\/ied classical
Yang--Baxter equation. In
Section~1 we recall the concept and basic facts of the $r$-matrix type Poisson
orbits. Then we describe the $r$-matrix Poisson pencil (i.e the pair of
compatible Poisson structures) of rank 1 or $CP^n$-type orbits of
$SL(n,C)$. Here we calculate symplectic leaves and the integrable foliation
associated with the pencil. We also describe the algebra of functions on
$CP^n$-type orbits. In Section~2 we calculate the Poisson homology of
Drinfeld--Sklyanin Poisson brackets which belong to the $r$-matrix Poisson
family.
\end{abstract}

\section*{Introduction}

The canonical or Poisson homology of Poisson manifolds were introduced by
Gelfand--Dorfman~\cite{kotov:GD}, Koszul~\cite{kotov:Kos} and Brylynsky~\cite{kotov:Br}.
Their algebraic
analogue was considered by Huebschmann in~\cite{kotov:Hue}.

Let $(M,\pi)$ be a smooth Poisson manifold with a Poisson structure given
by the antisymmetric bivector f\/ield $\pi$, and $f_i \in
  C^{\infty}(M)$, $i= 0,\ldots,k$ are smooth functions.

  Recall that the  formula for a canonical (Poisson) dif\/ferential of
  the degree $-1$ has  the following simple form on the decomposable
  dif\/ferential forms
\[
\ba{l}
\ds\delta_{\pi}(f_0 df_1\wedge \cdots \wedge df_k) =
\sum_i (-1)^{i+1} \pi(df_0,df_i)df_1\wedge \cdots \wedge{\hat {df_i}}\wedge
\cdots\wedge df_k
\vspace{3mm}\\
\ds \qquad
 +\sum_{i< j} (-1)^{i+j}  f_0 df_1\wedge \cdots
\wedge d\pi(df_i,df_j)\wedge \cdots \wedge
{\hat {df_j}}\wedge \cdots \wedge df_k.
\ea
\]

 The f\/irst def\/inition of \cite{kotov:GD} was inspired by the integrable systems
 theory: the Poisson homologies are responsable for the (non)-existence
 of bihamiltonian structures involved in the so-called Magri--Lenard
 scheme. It is now clear that there are many other reasons to study these
 homologies. We are unable to discuss all of them here and indicate only the
 most interesting points.

 The Poisson appeared in \cite{kotov:Br} as an important tool in the computations
 of Hochschild and cyclic homologies in the frame of Connes-like
 double complex, where the Poisson dif\/ferential $\delta_\pi$ plays the
 role of the Hochschild boundary operator and the usual de Rham
 dif\/ferential is very similar to the Connes cyclic cohomology
 operator.

 Moreover this ideology was used later by Feng--Tsygan~\cite{kotov:FT} for
 their description of the Hochschild complex
    for a ``quantum'' (deformed) algebra of smooth functions
     using the Poisson homology as the second term in the appropriate spectral
sequence.

  J.L. Brylinski, in the same paper \cite{kotov:Br}, conjectured an
  interesting symplectic version of the Hodge theorem. The negative
  answer to this conjecture~\cite{kotov:M,kotov:Yan} established a
  connection between the canonical cohomology complex of a
  symplectic manifold and the old topological problem,  where the
  homotopy type  of a space is a formal  consequence of its real
  homology ring. The interesting development of the conjecture and
  attempts to generalize it to the case of the Jacobi manifold, are
  given in the paper~\cite{kotov:FIdL}.

  O. Mathieu studied  various families of Poisson complexes~\cite{kotov:M1} and
  found that they are very dif\/ferent for dif\/ferent values of the
  parameter, even if the underlying
  Poisson structures are simply related (for example are parametrized
  by ${\mathbf CP_1}$). We would like to stress this aspect because of some
  similarities to our study of the Poisson complex associated to a
  Poisson pencil of $r$-matrix Poisson structures.

  O. Mathieu also linked the Poisson homology with the Gelfand-Fuks
  cohomology of the Lie algebra Hamiltonian vector f\/ields,
  symplectic operates and,
  f\/inally, his approach was used by G.~Papadopulo in his computations
  of cyclic (co)homology for Poisson and symplectic manifolds~\cite{kotov:P}.

  This list of interesting links would not be complete without mentioning
  the recent Weinstein def\/inition of the modular class of Poisson manifold,
  which can be considered as a classical analogue of the modular form
  of a von Neumann algebra~\cite{kotov:W,kotov:BrZ}. If this class is equal
  to $0$ the manifold is refered to as unimodular.

  The Poisson homology of an unimodular Poisson structure are
  isomorphic to the Poisson cohomology~\cite{kotov:Xu} and in fact there
  is an interesting pairing (which is a degenerate in general) between
  the Poisson homology and the Poisson cohomology~\cite{kotov:EvLuW}, whose
  algebraic roots are in the similar duality between the
Lie-Reinhart Poisson algebra cohomology and homology~\cite{kotov:Hue1}.

  We will f\/inish our brief survey of the recent manifestations of the
  Poisson homology with the indication of the link between the very
  general construction of Quillen which was used by B.~Fresse in his
  def\/inition of Poisson homology, and their computations for some
  singular Poisson surfaces~\cite{kotov:Fr}.

  In this text we consider the  Poisson algebraic structure associated
  with  a classical $r$-matrix, i.e. with a
  solution of the modif\/ied classical Yang--Baxter equation~\cite{kotov:D,kotov:STS}.
The classical $r$-matrix leads to Poisson orbits of two types.
Below a brief description for both of them is given.

The f\/irst structure, called the Drinfeld--Sklyanin, arises on Poisson
homogeneous spaces. It can be obtained as a result of Poisson reduction for
Poisson--Lie groups.

The second one exists on  special homogeneous spaces, which are known as
$r$-matrix orbits and classif\/ied in~\cite{kotov:GP}.
Hereafter  we will consider the  $r$-matrix type
orbits only.

It is clear that the Poisson brackets of $r$-matrix type arise both on the
Poisson--Lie groups and Poisson homogeneous spaces as the quasiclassical limit
  of the corresponding quantum objects.  Here we mean only a deformation
  quantization developed in~\cite{kotov:DGK,kotov:GRZ}.  Let us remark that, as
was shown in~\cite{kotov:K1}, the geometric quantization does not exist on certain
  Poisson homogeneous spaces.

   We now describe the structure of the paper.

  In Section 1 we recall the concept and basic facts of the $r$-matrix type
  Poisson orbits.  Then we describe the $r$-matrix Poisson pencil (i.e the pair
  of compatible Poisson structures) on the rank 1 or $CP^n$-type orbits of
  $SL(n,C)$. Here we calculate symplectic leaves and integrable foliation
  associated with the pencil. We also describe the algebra of functions on
  $CP^n$-type  orbits.

  In Section 2 we  calculate the Poisson homology of Drinfeld--Sklyanin
  Poisson brackets which belong to the $r$-matrix Poisson family.

  There are many interesting open questions the study of which we have
  postponed to the future. Among them the relations between our def\/inition
  of ``harmonic'' forms and the theory of Poisson harmonic forms, Kostant
  harmonic forms and equivariant Poisson cohomology~\cite{kotov:EvLu}, the
  precise links between Hochshild  and cyclic (co)homologies and the
  computations in this paper etc.

\section{On the Poisson structure of
             {\mathversion{bold}$CP^n$}-type complex orbits}

Let G be a semisimple Lie group and $\it g$
be a Lie algebra of G. Assume that $r$ is a standard Drinfeld--Jimbo $r$-matrix
\[
r=\sum_{\alpha \in \triangle^+} e_{\alpha}\wedge e_{-\alpha},
\]
where $e_{\alpha}$, $e_{-\alpha}$ is the Cartan basis for $\it g$.

Let $\cal O$ be a coadjoint orbit  and let $X_{\alpha}$,$X_{-\alpha}$
    be generators of the action of $G$
that corresponds
to the basis
$\{ e_{\alpha}, e_{-\alpha}\}$. Then the orbit $\mathcal O$ is called
     of $r$-matrix type if and only if the bivector f\/ield
\[
\pi=\sum_{\alpha \in \triangle^+} X_{\alpha}\wedge X_{-\alpha},
\]
 corresponds to the $r$-matrix $r$, gives us the
Poisson brackets.

If $G=SL(n)$ and $\mathcal O$ is an orbit of rank 1 matrices then the Poisson
 structure is called  $CP^n$-type.

Now let us give an explicit description of such orbits and $r$-matrix type
 brackets on them.

 Let us consider a standard action of $SL(n)$ on $C^n$ and its
cotangent lift to $T^{*}C^n$.  The
 generators of the action corresponding to the root basis $e_{ij}$ of $gl(n)$
 are in dual coordinates $(\zi,\xii)$
\[
X_{ij}=\zj\dzi-\xij\dxj
\]
The Drinfeld--Jimbo $r$-matrix
$r=\sum_{i<j}e_{ij}\wedge e_{ji}$ maps to the bivector f\/ield
\[
 \pi =\sum_{i<j}
 \zi \zj  \dzi \wedge \dzj - \sum_{i<j} \xii \xij  \dxi \wedge \dxj +\sum_{i<j} \zi
\xii \dzj \wedge \dxj -\sum_{j<i} \zi \xii  \dzj \wedge \dxj
\]
The structure is
 compatible with the natural symplectic one on $T^{*}C^n$ given by the form
\[
\omega = \sum_i d\zi \wedge d\xii.
\]

The momentum map $\mu:T^*C^n \rightarrow sl^*(n)$ is
\[
\mu(z,\xi)=\mbox{tr}\,(z^t\xi,*), \qquad \left(z^t\xi\right)_{ij}=\zi\xij .
\]
Thus the orbits of cotangent action of $SL(n)$ cover the orbits of rank 1 on
$sl^*(n)$. If $\mbox{tr}\,z^t\xi=(z,\xi)=\sum\limits_i \zi\xii \ne 0$
then the corresponding
orbit is semisimple and symmetric. If $\mbox{tr}\,z^t\xi=\sum\limits_i \zi\xii = 0$ and
$\mu(z,\xi) \ne 0$ then this is a nilpotent orbit of height 2.

The symplectic restriction to the level $(z,\xi)=\mbox{const}$ is a
pull-back of the Kirillov form. The momentum map is $SL(n)$-equivariant
so the $r$-matrix structure is compatible with the momentum map. Moreover the
momentum map gives us an isomorphism between  rank~1 orbits and  Poisson
reduction under the Hamiltonian action of $H=\sum\limits_i \zi\xii$.

Now we obtain the eigenvalues of the bivector $\pi_r$ with respect to the
symplectic form.
We introduce an operator f\/ield
$A$, $A(\phi)=V_{\Omega}(i_{\phi} \pi)$ for all 1-forms $\phi$. Here
$V_{\Omega}$ is a Hamiltonian operator  of the symplectic
structure $\Omega$ acting as
\[
V_{\Omega}: \dxi \rightarrow d\zi ,\ \dzi \rightarrow -d\xii
\]
First of all we check its eigenvalues and f\/ind it's eigen-vectors. Then we
``forget'' about the tangent direction to the f\/ibres of $\mu$, i.e. we separate
only the eigen-vectors that are the pull-back from the coadjoint orbits.

Let
\[
\phii =\xii d\zi +\zi d\xii, \qquad \psii =\xii d\zi - \zi d\xii .
\]

(i) We calculate how $A$ acts on the basis $\{\phii,\psii\}$ of 1-forms:
\[
\ba{l}
\ds i_{\phii}\pi =\sum_{i<j} (\xii\zi)\zj\dzj -\sum_{i<j} (\xii\zi)\xij\dxj
-\sum_{j<i} (\xii\zi)\zj\dzj +\sum_{j<i} (\xii\zi)\xij\dxj
\vspace{3mm}\\
\ds \qquad +\sum_{j<i} (\xij\zj)\xii\dxi -\sum_{j<i} (\xij\zj) \zi\dzi
-\sum_{i<j} (\xij\zj)\xii\dxi +\sum_{i<j} (\xij\zj)\zi\dzi,
\ea
\]
\[
\ba{l}
\ds i_{\psii}\pi =\sum_{i<j} (\xii\zi)\zj\dzj +\sum_{i<j} (\xii\zi)\xij\dxj
-\sum_{j<i} (\xii\zi)\zj\dzj -\sum_{j<i} (\xii\zi)\xij\dxj
\vspace{3mm}\\
\ds \qquad +\sum_{j<i} (\xij\zj)\xii\dxi +\sum_{j<i} (\xij\zj) \zi\dzi
-\sum_{i<j} (\xij\zj)\xii\dxi -\sum_{i<j} (\xij\zj)\zi\dzi.
\ea
\]
 Hence we see that
\[
\ba{l}
\ds A(\phii) =-\sum_{i<j} (\xii\zi)\zj d\xij -\sum_{i<j} (\xii\zi)\xij d\zj +
\sum_{j<i} (\xii\zi)\zj d\xij +\sum_{j<i} (\xii\zi)\xij d\zj
\vspace{3mm}\\
\ds \qquad +\sum_{j<i} (\xij\zj)\xii d\zi +\sum_{j<i} (\xij\zj) \zi d\xii
-\sum_{i<j} (\xij\zj)\xii d\zi -\sum_{i<j} (\xij\zj)\zi d\xii.
\ea
\]

If we def\/ine $a_i =\xii \zi$ then
\[
A(\phii) =\left(\sum_{j<i} a _j -\sum_{i<j} a_j\right)\phii +
a_i \left(\sum_{j<i} \phi_j -\sum_{i<j} \phi_j\right) .
\]
Similarly
\[
A(\psii) =\left(\sum_{j<i} a_j -\sum_{i<j} a_j\right)\psii
-a_i \left(\sum_{j<i} \psi_j -\sum_{i<j} \psi_j\right).
\]

(ii) Note that if $\phi_0 =\sum\limits_i \phii$, then
\[
A(\phi_0)=\sum_i\left(\sum_{j<i} a_j -\sum_{i<j} a_j\right)\phii
+\sum_i a_i \left(\sum_{j<i} \phi_j -\sum_{i<j} \phi_j\right)=0.
\]

\begin{explanation}
This form vanishes on the level $\sum\limits_i \zi\xii=\mbox{\rm const}$.
Therefore, if we make a
Poisson reduction via the Hamiltonian field $H=\sum\limits_i \zi\xii$,
then one needs to consider $A\; \mbox{\rm mod}\,(\phi_0)$.
\end{explanation}
Thus
\[
\ba{l}
\ds  A(\phii)=\left(\sum_{j< i} a_j -\sum_{i<j} a_j\right)\phii
+ a_i \left(\sum_{j<i} \phi_j -\sum_{i<j}\phi_j\right)
\vspace{3mm}\\
\ds \qquad = \left(\sum_{j\le i} a_j -\sum_{i<j} a_j\right) \phii
+ a_i \left(-\phi_0 +2\sum_{j<i} \phi_j\right).
\ea
\]
We obtain the ``triangle'' basis $\{ \phi_0 ,\phi_1,\ldots,\phi_{n-1},
 \psi_0 ,\psi_1,\ldots,\psi_{n-1}\}$ for $A$:
\[
\ba{l}
\ds A(\phi_0)=0; \\
\cdots \cdots \cdots\cdots\cdots \\
\ds A(\phii)=\left(\sum_{j\le i} a_j -\sum_{i<j}
a_j\right) \phii+ a_i \left(-\phi_0 +2\sum_{j<i} \phi_j\right);
\vspace{3mm}\\
\cdots \cdots \cdots\cdots\cdots\\
\ds A(\phi_{n-1})=\left(\sum_{j< n} a_j -a_n\right) \phi_{n-1}
+ a_{n-1} \left(-\phi_0 +2\sum_{j<n-1} \phi_j\right).
\ea
\]

We see that the eigenvalues of $A$ are  equal to
$\sum\limits_{j\le i} a_j -\sum\limits_{i<j}a_j$, $i=1,\ldots,n-1$ or
$\lambda_i=\sum\limits_{j\le i}\zj\xij -\sum\limits_{i<j} \zj\xij$,
$i=1,\ldots,n-1$.
Moreover,
\[
\ba{l}
\ds A\left(\sum_{i\le k}\phi_i\right), \mbox{mod}(\phi_0)=\sum_{i\le k}
\sum_{j\le i} a_j\phi_i-\sum_{i\le k}\sum_{i<j} a_j\phi_i +
2\sum_{i\le k}\sum_{j<i} a_i\phi_j
\vspace{3mm}\\
\ds \qquad =\sum_{i\le k}\sum_{j\le i} a_j\phi_i-
\sum_{i< k}\sum_{j\le k} a_j\phi_i-
\sum_{i\le k}\sum_{j >k} a_j\phi_i+
2\sum_{j\le k}\sum_{i<j} a_j\phi_i=
\vspace{3mm}\\
\ds \qquad =\sum_{i\le k}\sum_{j\le i} a_j\phi_i+
\sum_{j\le k}\sum_{i<j} a_j\phi_i-
\sum_{i\le k}\sum_{j >k} a_j\phi_i=
\sum_{i\le k}\sum_{j\le k} a_j\phi_i-
\sum_{i\le k}\sum_{j >k} a_j\phi_i=
\vspace{3mm}\\
\ds \qquad =\left(\sum_{j\le k}a_j -\sum_{j>k}a_j\right)
\sum_{i\le k}\phi_i.
\ea
\]

Hence the diagonal $\mbox{mod}(\phi_0)$-basis is
\[
\bph_k=\sum_{i\le k}\phi_i, \qquad
A(\bph_k)=\lambda_k\bph_k, \qquad  k=1,...,n-1.
\]
 Note that $\bph_k=d\left(\sum\limits_{i\le k}\zi\xii\right)$.

It is not dif\/f\/icult to see also that
$A(\bps_k)=\lambda_k\bps_k$, $k=1,...n-1$, where
\[
\bps_k=\frac{\psi_{k+1}}{a_{k+1}}-
\frac{\psi_{k}}{a_{k}}.
\]
 Note that
\[
\bps_k=d(\ln(z_{k+1})-\ln(\xi_{k+1})-\ln(z_k)+\ln(\xi_k))=
d\left(\ln\left(\frac{\xi_k z_{k+1}}{z_k\xi_{k+1}}\right)\right).
\]
\begin{explanation}
All forms $\bph_k$, $\bps_k$  vanish on the orbits of the $H$-Hamiltonian
action
and are $H$-invariant, so one has to consider them as a pull-back of
the almost
everywhere independent basis of differential 1-forms from coadjoint
orbits.
\end{explanation}

Hence we get an integrable foliation on the orbits associated with
a Poisson pencil $\pi_r +\lambda\pi_{kir}$,
where $\pi_{kir}$ is the Kirillov Poisson structure that is also called
``Lie--Poisson''.
 It is well-known that $A$, as def\/ined above, is ``integrable''
since $\pi_r$ and $\pi_{kir}$ are compatible structures. This means that  its
eigen- (or adjacent) spaces give us an integrable (singular) foliation.  As
proved above, the foliation is given by the equations $\bph_k=0$, $\bps_k=0$,
$k\in I\subset (1,\ldots,n-1)$ or
\[
 \sum_{i\le k}\zi\xii=c_k,\qquad  \frac{\xi_k
z_{k+1}}{z_k\xi_{k+1}}=c'_k, \quad  k\in I.
\]


We know at least two natural Poisson algebras associated with
an $r$-matrix on the rank~1 orbits.

 The f\/irst one is an algebra $A$ generated by the Hopf bundle on $CP^n$.
It might be identif\/ied with the algebra of polynomials
$C[z_1,\ldots,z_n]$, where $SL_n$
acts as follows:
\[
e_i=z_i\partial_{z_{i+1}}, \qquad
 f_i=z_{i+1}\partial_{z_i}, \qquad  h_i=z_i\partial_{z_i}-z_{i+1}
\partial_{z_{i+1}}.
\]
 Here $e_i$, $f_i$, $h_i$  are Cartan generators of the Lie algebra.

 Recall that
\[
 C^k[z_1,\ldots, z_n] \cong V_{\omega_1},
\]
 where $C^k[z_1,\ldots, z_n]$ is the space of polynomials of degree $k$,
   $V_{\omega_1}$ is an $sl_n$-module with highest weight $\omega_1$.

 The second one is the algebra $B$ of algebraic (holomorphic)
functions on $\mathcal O$.
  The algebra of functions on $CP^n$-type orbits
is generated by polynomials $\zi\xij$ with the relation
$\sum\limits_i \zi\xii =\mbox{const}\ne 0$.

The main facts concerning the structure of the algebra of functions
on $\mathcal O$ as
$sl_n$-module are as follows:
\[
\ba{ll}
1. &  \mbox{Fun}({\mathcal O})\cong {\mathop{\oplus}\limits_{k\ge 0}}
 V_{k(\omega_1 + \omega_{n-1})},
\vspace{2mm}\\
2. & C^k[z]\otimes C^k[\xi] \stackrel{\mathrm{def}}{=} C^{k,k} [z,\xi]
\cong {\mathop{\oplus}\limits_{l\le k}} V_{l(\omega_1 + \omega_{n-1})}.
\ea
\]

\begin{proposition}The subcomplex of algebraic differential forms
$\Omega^k ({\mathcal O})$
on ${\mathcal O}$, generated by  subalgebra  $B^k({\mathcal O})\cong
{\mathop{\oplus}\limits_{l\le k}} V_{l(\omega_1 + \omega_{n-1})} $
 of $\mbox{\rm Fun}({\mathcal O})$,  is isomorphic to the subcomplex
 $A^k$  of forms  on $T^*C^n$, generated by $\zi\xij$.
\end{proposition}

\noindent
{\bf Proof.}
As we know from the Hochschild--Konstant--Rosenberg theorem \cite{kotov:HKR}
\[
\Omega^*(B)\cong HH^*(B),
\]
 where $HH^*(B)=H^*(C^*(B),b)$
are Hochschild homologies of $B$. They are calculated as the homology of
complex of
$B$-chains $C^*(B)=\oplus_m B\otimes B^{\otimes m}$, with dif\/ferential $b$
acting as follows
\[
\ba{l}
\ds b(a_0\otimes a_1\ldots \otimes a_m)
\vspace{3mm}\\
\ds \qquad =\sum_{j\le m-1}
(-1)^j a_0\otimes \ldots a_j a_{j+1} \ldots \otimes a_m+
(-1)^m a_m a_0\otimes a_1\ldots \otimes a_{m-1}.
\ea
\]

Poisson and Hochschild dif\/ferentials on $B$-chains
are compatible with the introduced f\/iltration, so we can deduce that
\[
A^k \cong HH^{k,k}_*(C[z,\xi])\cong HH^k_*(B)\cong
\Omega^k(\mathcal O).
\]

Let $\mathcal O$ be a $CP^n$-type orbit and $P={\mathop{\oplus}\limits_k} P^k$
 -- subalgebra of
polynomials on $\zi$, $\xij$, generated by $\zi \xij$. Let $B$ be an algebra
of functions on $\mathcal O$. Then, without losing generality, we establish
that
\[
B=P/J, \qquad  J=\{ p(\xi, z)=(H(\xi, z)-1) f(\xi, z)| f\in P \},
\]
where   $H(\xi, z)=\sum\limits_i \zi\xii$.

\newpage

 The algebra $B$ has a natural f\/iltration $B^{(k)}$, arising from the
f\/iltration on $P$, $P^{(k)}={\mathop{\oplus}\limits_{j\le k}} P^j$.
\[
 B^{(k)} = P^{(k)}/ (H-1)P^{(k-1)}.
\]

\begin{proposition} We have the following isomorphism:
\[
B^{(k)}\cong P^k.
\]
\end{proposition}

\noindent
{\bf Proof.}
Let $f\in B^{(k)}$ and $f=\sum\limits_{j=0}^{k} f_j$
be a decomposition on homogeneous components.
Consider the mapping
\[
{\mathcal P}: f \mapsto  \hat f =\sum f_j H^{k-j}, \  B^{(k)}\rightarrow P^k.
\]
Then $\mbox{Ker}\; {\mathcal P}=J^{(k)}=(H-1)B^{(k-1)}$. Indeed
\[
 \hat f =\sum f_j H^{k-j}=f+\sum_{j<k} (H^{k-j}-1)f_j =0 \
\Longrightarrow \  f\in J^{(k)}
\]
so $\mbox{Ker}\; {\mathcal P} \subset J^{(k)}$.
On the other hand, if $f\in J^{(k)}$ then
\[
f_k=Hg_k, \qquad f_j=Hg_{j-1}-g_j.
\]
Therefore it is easy to see that $\mbox{Ker}\; {\mathcal P} = J^{(k)}$.

Thus ${\mathcal P}$ gives rise to the isomorphism of f\/iltered spaces
\[
{\cal P}: B^{(k)}\rightarrow P^k .
\]

\begin{proposition}
The map ${\mathcal P}$ satisfies the following conditions:

1. ${\mathcal P}$ is a homomorphism of algebras, i.e.
\[
\hat {f_1 f_2} =\hat f_1 \hat f_2;
\]

2.  ${\mathcal P}$ is a homomorphism of $SL_n$-modules, i.e.
 the following diagramm is commutative
\[
g\in SL_n, \qquad
\begin{array}{ccc}
B^{(k)} & \rightarrow & P^k \\
\downarrow g && \downarrow g \\
B^{(k+1)} & \rightarrow & P^{k+1}.
\end{array}
\]
\end{proposition}

\noindent
{\bf Proof.}
1. Let $f_1=\sum\limits_{j=0}^k f_{1,j}$,
$f_2=\sum\limits_{j=0}^l f_{2,j}$. Then
\[
\hat f_1=\sum_{j=0}^k f_{1,j} H^{k-j} ,\qquad \hat f_2=\sum_{j=0}^l f_{2,j}
H^{l-j}.
\]
Hence
\[
\hat {f_1 f_2} =\sum_{i+j\le k+l} f_{1,i} f_{2,j} H^{k+l-i-j}=\hat f_1 \hat
f_2.
\]

2. It follows from the fact that $H$ is invariant under the $SL_n$-action.

 It is very important that two  algebras of complex smooth functions on
$CP^{n}$ and of algebraic functions on $\mathcal O$ are isomorphic
as $sl_n$-modules. So the following diagrams are commutative $$
\begin{array}{ccc}
\Omega^*(CP^{n-1})&\rightarrow & \Omega^*({\mathcal O})
\vspace{1mm}\\
\downarrow d&& \downarrow d
\vspace{1mm}\\
\Omega^{*+1}(CP^{n-1})&\rightarrow & \Omega^{*+1}({\mathcal O})
\end{array}
$$
$$
\begin{array}{ccc}
\Omega^*(CP^{n-1})&\rightarrow & \Omega^*({\mathcal O}) \vspace{1mm}\\
\downarrow \delta_{\omega},\dep && \downarrow \delta_{\omega},\dep
\vspace{1mm}\\
\Omega^{*-1}(CP^{n-1})&\rightarrow & \Omega^{*-1}({\mathcal O}),
\end{array}
$$
where $d$ is the de Rham dif\/ferential and $\dep$, $\delta_{\omega}$ are Poisson
dif\/ferentials, corresponding to the $r$-matrix structure $\pi$ and to the
symplectic form $\omega$, respectively.

As a direct consequence we get the identity
\[
 H_{*,\dep'}(CP^n)\cong H_{*,\dep '} ({\mathcal O})
\]
for all Poisson dif\/ferentials $\dep '=a\delta_{\omega}+b\dep$
from the Poisson pencil $\pi$, $\pi_{\omega}$.

We remark that the structure of the $sl_n$-module on the algebra of complex
functions arises as a simple complexif\/ication of the
$su_n$-module structure on the smooth function algebra.

Let
\[
\{ f_1, f_2 \}\stackrel{\mathrm{def}}{=} \langle \pi, df_1\wedge df_2\rangle,
\qquad
\{ f_1, f_2 \}_{\omega}\stackrel{\mathrm{def}}{=} H \langle \pi_{\omega}, df_1\wedge df_2
\rangle
\]
and $\{f_1, f_2\}_{a,b}=a\{f_1, f_2\}+b\{f_1, f_2\}_{\omega}$. Then we can
prove the  following proposition.

\begin{proposition}
The following identity holds:
\[
\hat{\{ f_1, f_2\}}_{a,b} = \{\hat f_1,\hat f_2\}_{a,b}.
\]
\end{proposition}

\noindent
{\bf Proof.}
Straightforward calculations show that
the mapping ${\mathcal P}$
satisf\/ies this formula because $H(z,\xi)$ is Casimir
with respect to the whole pencil $\pi_{a,b}$.
The structure $\pi_{\omega}$ is symplectic and nongenerate on $T^* C^n =\{ (z,
\xi \}$. We easily obtain $\{H, \zi\xij\}_{\omega}=0$.


 \section{On the Poisson homology of
          {\mathversion{bold}$CP^n$}-type  orbits}

The Poisson homology was introduced as the second
term in the spectral sequence associated with the Hochschild complex for a
deformed algebra of smooth functions.

We give a short description of the Poisson homology.

Let $X$ be a smooth manifold and $A_0=C^{\infty} (X)$ be an algebra of smooth
functions on $X$. An associative algebra $A$ over the ring of formal series
$C[[h]]$ is called a deformation of $A_0$ if $A$ is isomorphic to $A_0\otimes
C[[h]]$ as a $C[[h]]$-module and a multiplication on $A$ coincide with the
multiplication on $A_0\otimes C[[h]] \;\mbox{mod}\, O(h)$~\cite{kotov:Lch}.

Every such deformation can be obtained by deformation
quantization~\cite{kotov:Lch} of Poisson brackets  on $X$, which means that
the commutator on the deformed algebra give us a Lie algebra structure
on the functions on~$X$
\[
(a*b-b*a) \;\mbox{mod} \,o(h) =\{a,b\}.
\]

Recently it was shown by M.~Kontsevich~\cite{kotov:Ko}
that every Poisson structure on a f\/lat space is quantizible.

It is well-known fact  that there exists a
natural chain complex associated with any
associative algebra. Here  we recall the def\/inition and some basic facts about
this complex.

Let  $A$ be an associative algebra with unity over the f\/ield $C$.

Let us denote $C_n(A)=\cn$ as the space of $A$-value
$n$-chains. We can consider it
as $A$-bimodule with the usual left and right actions
\[
\ba{rcl}
( a,{a_0\otimes a_1\otimes \ldots \otimes a_{n-1}\otimes a_n} )&\rightarrow & a
{a_0\otimes a_1\otimes \ldots \otimes a_{n-1}\otimes a_n} \vspace{2mm}\\
(a,{a_0\otimes a_1\otimes \ldots \otimes a_{n-1}\otimes a_n})
 &\rightarrow &{a_0\otimes a_1\otimes \ldots \otimes a_{n-1}\otimes a_n} a.
\ea
\]

 Recall that there
is an  exact sequence of right $A$-modules
\[
C_{n+1}(A) \stackrel{b'}{\rightarrow} C_n(A) \stackrel{b'}{\rightarrow}
 \cdots \stackrel{b'}{\rightarrow} C_1(A) \stackrel{b'}{\rightarrow} A,
\]
where $b'=\sum\limits_{0\le i\le n-1} b_i$ and
\[
b_i ({a_0\otimes a_1\otimes \ldots \otimes a_{n-1}\otimes a_n})=
(-1)^i a_0\otimes \cdots \otimes a_i a_{i+1}\otimes \cdots \otimes a_n.
\]
 We introduce an  operator of homotopy $s:C_{*}\rightarrow C_{* +1}$,
\[
 s({a_0\otimes a_1\otimes \ldots \otimes a_{n-1}\otimes a_n})=1\otimes
{a_0\otimes a_1\otimes \ldots \otimes a_{n-1}\otimes a_n},
\]
such that the identity $b's+sb'=1$ holds. So the complex $(C(A),b')$ is
acyclic.

We say that an associative algebra is $H$-unital if the complex $(C(A),b')$
is homotopicaly trivial.

 On can consider another operator $b'=\sum\limits_{0\le i\le n} b_i$, where
 $b_n ({a_0\otimes a_1\otimes \ldots \otimes a_{n-1}\otimes a_n})=
(-1)^n a_n a_0\otimes a_1\otimes \cdots  \otimes a_{n-1}. $

 We see again that $b^2=0$, so there is a new complex named the Hochschild
ones
\[
C_{n+1}(A) \stackrel{b}{\rightarrow} C_n(A) \stackrel{b}{\rightarrow}
 \cdots \stackrel{b}{\rightarrow} C_1(A) \stackrel{b}{\rightarrow} A.
\]

 The homologies of this complex are the Hochschild homologies.

One can also def\/ine it as $\mbox{Tor}_{A\otimes A^o} (A,A)$ in the category of
$A$-bimodules or $A\otimes A^o$-modules, where $A^o$ is an ``opposite'' algebra.
 We used the free resolution def\/ined above to establish an
explicit formula for the Hochschild complex. The result of course does not
depend on the choice of the projective resolution.

The dual construction $HH^n (A,A)=\mbox{Ext}^n_{A\otimes A^o} (A,A)$ is called a
Hochschild cohomology.

\medskip

\noindent
{\bf Example 1} (Hochschild--Kostant--Rosenberg~\cite{kotov:HKR})
Let $A=C^{\infty} (X)$, where $X$ is a smooth manifold. Then $HH_n (A,A)\simeq \Omega_n
(X)$.
Two maps
\[
\chi:{a_0\otimes a_1\otimes \ldots \otimes a_{n-1}\otimes a_n}
\rightarrow \frac{1}{n!} a_o da_1\wedge \cdots\wedge da_n
\]
and
\[
\chi^{-1} : a_o da_1\wedge \cdots\wedge da_n \rightarrow
a_0\otimes\sum_{\sigma\in S_n}\epsilon(\sigma) a_{\sigma1}\otimes \cdots
\otimes a_{\sigma n}
\]
give us a quasi-isomorphism of complexes.

\medskip

Now let $A$ be the deformed algebra of smooth functions on $X$ introduced
above. Let us consider a f\/iltration for the Hochschild complex of $A$ def\/ined
as follows:
\[
F^p C_*(A)=h^p C_*(A)
\]
Note that the f\/iltration is served by the Hochschild dif\/ferential and there
is an isomorphism  $F^p C_*(A)/F^{p+1} C_*(A)\simeq A_0$ as vector spaces.

Apparently the zero-order term in the spectral sequence associated with the
f\/iltration coincides with the Hochschild complex for $A_0$ as a
free $C[[h]]$-module. That is why we get (from the  Hochschild--Kostant--Rosenberg
theorem)
\[
E_0^p (C_*(A))\simeq \Omega^* (X) \otimes C[[h]].
\]

As shown in \cite{kotov:Br}, the f\/irst-order term in the corresponding
sequence coincides with the Poisson homology complex $(\Omega^*(X),
\delta_{\pi})\otimes C[[h]]$, where the dif\/ferential $\dpi$ is def\/ined as
follows:
\[
\ba{l} \dpi(f_0 df_1\wedge \cdots \wedge df_k)=\sum_i (-1)^{i+1}
\{f_0,f_i\} df_1\wedge \cdots \wedge  \hat{df_i}\wedge
\cdots \wedge df_k
\vspace{2mm}\\
\ds \qquad +\sum_{i<j}(-1)^{i+j}f_0
df_1\wedge \cdots \wedge \hat{df_i}\wedge \cdots\wedge \hat{df_j}\wedge \cdots
\wedge df_k
\ea
\]
or, equivalently
\[
 \dpi(\tau)=di_{\pi}(\tau)-i_{\pi}d(\tau),\qquad  \tau\in\Omega^*.
\]


Now the task is to study the Poisson homology of $r$-matrix type
coadjoint orbits of rank~1.
The programme of investigation for  the Hochschild and
cyclic homology of the $CP^n$-type orbits is supposed.

Recall that
if $\pi=\sum\limits_i X_i\wedge Y_i$ then
$\dpi=\sum\limits_i (L_{X_i}i_{Y_i}-i_{X_i}L_{Y_i})$,
where  $L_X$ and $i_X$ are the Lie derivatives along $X$ and
the interior multiplication by $X$, respectively.


\medskip

\noindent
{\bf Example 2}
 The $r$-matrix structure introduced above might be written on
the  $CP^n$-type
orbits as a restriction of the bivector f\/ield
\[
 \pi =-\sum_{i<j} \zi \zj  \dzi \wedge \dzj + \sum_{i<j} \xii \xij  \dxi \wedge
\dxj
-\sum_{i<j} \zi \xii \dzj \wedge \dxj +\sum_{j<i} \zi \xii  \dzj \wedge
\dxj
\]
to the space of forms on $(z,\xi)$ of equal degree with the only relation
$\langle z,\xi\rangle =\sum\limits_i \zi\xii =\mbox{const}$. We assume
that
$\langle z,\xi\rangle =\sum\limits_i \zi\xii =1$.
Hence the Kirillov--Kostant--Souriau structure comes from the restriction
of the bivector f\/ield
\[
 \pi_{\omega}=\langle z,\xi\rangle \sum_i \dzi\wedge\dxi.
\]

That is why
\[
\ba{l}
\ds \pi_{a,b}=a\pi_{\omega}+b\pi=
a\sum_i \zi\xii \dzi\wedge\dxi -b
\sum_{i<j} \zi \zj  \dzi \wedge \dzj + b\sum_{i<j} \xii \xij  \dxi \wedge
\dxj
\vspace{3mm}\\
\ds \qquad
+(a+b)\sum_{i>j} \zi \xii \dzj \wedge \dxj +(a-b)\sum_{i<j} \zi \xii  \dzj \wedge
\dxj.
\ea
\]

Now we compute a Poisson homology for
\[
\ba{l}
\ds \pi_{DS}=\pi_{\omega}+\pi=
\sum_i \zi\xii \dzi\wedge\dxi
\vspace{3mm}\\
\ds \qquad -
\sum_{i<j} \zi \zj  \dzi \wedge \dzj + \sum_{i<j} \xii \xij  \dxi \wedge
\dxj +2\sum_{i>j} \zi \xii \dzj \wedge \dxj.
\ea
\]

One can introduce a grading on the algebras of algebraic forms and polyvectors
with formal coef\/f\/icients such that all natural operations (i.e
exterior dif\/ferential and multiplication, interior multiplication and
Schouten--Nijenhuis brackets) preserve the grading. Here
we attach the grading degree $i$ to
$\zi$, $\xii$, $d\zi$, $d\xii$ and $-i$ to $\dzi$, $\dxi.$

The algebra of formal dif\/ferential forms
$\Omega^*_{\mbox{\scriptsize form}} (\zi,\xii)$ is graded (as
well as the algebra~$A^k$). The corresponding increasing f\/iltration
is
\[
 F^p A=\bigoplus_{k\ge p}A_p, \qquad A=F^0 A\supset F^1 A\supset \cdots \supset
F^p A\supset \cdots .
\]

The Poisson (Brylinsky) dif\/ferential $\dpi$  agrees with this
f\/iltration. In order to see it we decompose $\pi_{DS}$ in two
components, $\pi_{DS}=\pi_0+\pi_1$, where
\[
\ba{l}
\ds \pi_0=
-\sum_{i<j} \zi\zj\dzi\wedge\dzj + \sum_{i<j} \xii\xij\dxi\wedge\dxj +
\sum_i \zi\xii \dzi\wedge\dxi
 \vspace{3mm}\nonumber\\
\ds \pi_1= 2\sum_{i>j} \zi \xii \dzj \wedge \dxj.
\ea\nonumber
\]

The f\/irst component of $\pi_{DS}$ is of degree $0$, the second is the sum
of bivector f\/ields of positive degrees. So $\dpi (F^p A)\subset F^p A$.
Moreover in the corresponding spectral sequence~$E^*_r$
\[
E^*_0=\bigoplus_p F^pA / F^{p+1}A\simeq\bigoplus_p A_p= A, \qquad
 \delta_o= \delta_{\pi_0}. 
\]

 One has to compute the Poisson homology for $\pi_0$.
The Poisson dif\/ferential $\delta_{\pi_0}$ is
\[
\ba{l}
\ds \delta_{\pi_0}=
\sum_{i<j}(-\Lzi\izj+\izi\Lzj+\Lxi\ixj-\ixi\Lxj)
\vspace{3mm}\\
\ds \qquad +\sum_{i}(\Lzi\ixi-
\Lxi\izi)
=\sum_j\left(\!-\sum_{i<j}\Lzi+\sum_{j<i}\Lzi -
\Lxi\!\right)\izj
\vspace{3mm}\\
\ds \qquad -\sum_j \left(-\sum_{i<j}\Lxi+\sum_{j<i}\Lxi -\Lzj\right)\ixj.
\ea
\]

Now our task is to compute the Poisson homology for the subcomplex of forms
of equal total degree on $z$ and $\xi$ that in addition are in the kernel of
$i_H \stackrel{\mathrm{def}}{=} i_{X_H}$, where $X_H=\sum\limits_i
\zi\dzi -\sum\limits_i\xii\dxi$.

Let us denote $A^k (m,l)=\{ \omega\in A^k | \; \mbox{deg}_{\zi} \omega =m_i,
\; \mbox{deg}_{\xii} \omega =l_i     \}$, $m, l \in Z_+^n$.

It is clear that $A^k (m,l)$ is a subcomplex of $A^k$ under
$\delta_{\pi_0}$.  On the subcomplex this dif\/ferential can be written as
\[
\delta_{\pi_0}
=\sum_j\left(-\sum_{i<j}m_i+\sum_{j<i}m_i -l_j\right)\izj
-\sum_j \left(-\sum_{i<j}l_i+\sum_{j<i}l_i -m_j\right)\ixj,
\]
or
\[
\ba{l}
\ds \delta_{\pi_0} =\sum_j a_j\izj
-\sum_j b_j\ixj,
\vspace{3mm}\\
\ds a_j=-\sum_{i<j}m_i+\sum_{j<i}m_i -l_j, \qquad
b_j=-\sum_{i<j}l_i+\sum_{j<i}l_i -m_j.
\ea
\]

Now we f\/ind the homotopy operators in the form
\[
s=\sum_j x_j\zj^* d\zj + \sum_j p_j\xij^* d\xij  =\sum_j x_j\izj^* +
\sum_j p_j\ixj^*
\]
such that
\[
si_H +i_H s=0.
\]

In this case $s: \mbox{Ker}\;  i_H \rightarrow \mbox{Ker} \; i_H$.

>From the well known commutation relations
\[
\ba{lll}
\{d\zi, d\zj^* \}=d\zi d\zj^* +d\zj^* d\zi =\delta_{i,j}, &\quad &
 [\zi^*, \zj ]=\zi^*\zj -\zj\zi^* =\delta_{i,j},
\vspace{2mm}\\
\{d\xii, d\xij^* \}= d\xii d\xij^* +d\xij^* d\xii=\delta_{i,j} & &
  [\xii^*, \xij ]=\xii^* \xij -\xij \xii^* =\delta_{i,j}
\ea
\]
 we obtain that
\[
\{ \izi^*,\izj\}=\Lzi\delta_{i,j}= m_i\delta_{i,j}, \qquad
 \{ \ixi^*,\ixj\}=\Lxi\delta_{i,j}=  l_i\delta_{i,j}
\]
and
\[
\{ s, i_H\}=si_H +i_H s=0 \ \Longleftrightarrow \ \sum_i x_i m_i -\sum_i p_i l_i=0.
\]

>From the identity
\[
\delta_{\pi_0}s+s\delta_{\pi_0}=\sum_i a_j x_j m_j -\sum_i b_j p_j l_j
\]
we see that the subcomplex $A^k(m,l)$ is acyclic if there
exist such $(x_1,\ldots, x_n, p_1, \ldots, p_n)$ that
\[
\sum_i a_j x_j m_j -\sum_i b_j p_j l_j=0, \qquad
\sum_i x_i m_i -\sum_i p_i l_i=0.
\]

This is false  if and only if  for some $\lambda$
\[
a_j m_j =\lambda m_j, \qquad b_j l_j =\lambda l_j.
\]
 Let
\[
I=\{ i | m_i \ne 0 \}, \qquad J=\{j | l_j \ne 0 \}
\]
Thus we have the equation
\[
\left\{
\ba{ll}
m_i=l_j=0, & i\in I, \ j\in J
\vspace{1mm}\\
a_i=b_j, & i\in \{ 0,\ldots, n\}-I, \ j\in
\{ 0,\ldots, n\}-J
\ea
\right.
\]
or
\[
\left\{
\ba{l}
m_i=l_j=0,  \qquad  i\in I, \ j\in J
\vspace{3mm}\\
\ds \sum_{i<j}m_i-\sum_{j<i}m_i +l_j =
\sum_{i<k}l_i-\sum_{k<i}l_i +m_k,
\vspace{3mm}\\
\ds  i\in \{ 0,\ldots, n\}-I,\quad k\in \{ 0,\ldots, n\}-J.
\ea
\right.
\]

Now we can  solve it.
Let $ m_j, m_k \ne 0$, $j<k$. Then $a_j=a_k$ and
\[
\ba{l}
\ds \sum_{i<j}m_i-\sum_{j<i}m_i +l_j =
\sum_{i<k}m_i-\sum_{k<i}m_i +l_k \quad \Longrightarrow
\vspace{3mm}\\
\ds \hspace*{6cm}
\Longrightarrow \quad 2\sum_{j<i<k} m_i +m_j +m_k -l_j +l_k =0.
\ea
\]

Thus $l_j>l_k\ge 0$.

Let $l_k>0$. Then $a_j=a_k=b_j=b_k$ and
\[
\left\{
\ba{l}
\ds 2\sum_{j<i<k} m_i +m_j +m_k +l_j -l_k =0
\vspace{3mm}\\
\ds 2\sum_{j<i<k} l_i -m_j +m_k +l_j +l_k =0
\ea
\right. \Longrightarrow\
2\sum_{j<i\le k} (m_i +l_i) =0 \ \Longrightarrow \ l_k=m_k=0.
\]

Therefore, $l_k=0$. That's why $I\cap J=\{j\}$, $\# (I\cap J)=1$.

Let $ i<j<k$, $m_i, m_j, m_k \ne 0$. Then $l_i>l_j>l_k\ge 0$ and
$\#(I\cap J)\ge 2$. So there exists no more than $2$
dif\/ferent indexes $i$,$j$
such that $m_i, m_j >0$. Hence  $\#I, \#J \le 2$.

Let $m_i, m_j, l_j, \ne 0$. Then we obtain $a_i=a_j=b_j$
\[
\mbox{if} \ i<j \quad
\left\{
\ba{l}
a_i=m_j\\
a_j=-m_i-l_j
\ea
\right. \
a_i=a_j=b_j \quad  \mbox{then}  \quad m_i+m_j+l_j=0
\]
and
\[
\mbox{if}\  i>j, \quad
\left\{
\ba{l}
a_i=-m_j\\
a_j=m_i-l_j \\
\ea
\right. \
a_i=a_j \quad \mbox{then}  \quad m_i+m_j=l_j.
\]

Since $\sum m_i=\sum l_i$ we get $l_i=0$, $i\ne j$, and
the solution is $I=\{j,i |j<i \}$, $J=\{j \}$ and
$J=\{j,i |j<i \}$, $I=\{j \}$.

In this case
\[
\delta_{\pi_0}=i_H=\sum_i \izi -\sum_i \ixi =0
\]
on the subcomplex $A^k(m,l)$.

Now one can declare that
\[
\ba{l}
\ds H^*_{\delta_{\pi_0}} \cong \lim_{k\rightarrow \infty} \Bigl\{ \omega \in A^k|\omega
=\sum_{i\ge j} \sum_{0\le\mu \le 1, 0\le\nu\le 1}
\sum_{p+q+\mu+\nu}
\alpha_{i,j,\mu,\nu} (\zi\xij)^p d(\zj\xij)^{\mu} \}
\vspace{3mm}\\
\ds \hspace*{8cm} +\beta_{i,j,\mu,\nu} (\zj\xii)^q d(\zj\xij)^{\nu} \Bigr\}.
\ea
\]

Hence one can calculate the f\/irst term in the spectral sequence.
Now let us prove that the spectral sequence converges in the f\/irst term.

Since $(\zi \xij)^p d(\zi\xij)^{\mu} (\zj \xij)^q d(\zj\xij)^{\nu}$, $i>j$ is
of weight $i(p+\mu)+j(p+\mu+2q+2\nu)$, it is suf\/f\/icient to f\/ind such cocycle
$\sigma$ that
\[
(\zi \xij)^p d(\zi\xij)^{\mu} (\zj \xij)^q d(\zj\xij)^{\nu}-\sigma=0
\quad \left(\mbox{mod}\;
F^{i(p+\mu)+j(p+\mu+2q+2\nu)+1} \right).
\]

Let
\[
\sigma=(\zi \xij)^p d(\zi\xij)^{\mu} \left(\sum_{k\ge j}z_k \xi_k\right)^q
d\left(\sum_{k\ge j}z_k\xi_k\right)^{\nu}.
\]

Then
\[
\ba{ll}
1. & \delta_{\pi}\sigma=0;
\vspace{2mm}\\
2. & \ds \sigma \in F^{i(p+\mu)+j(p+\mu+2q+2\nu)};
\vspace{2mm}\\
3. &
(\zi \xij)^p d(\zi\xij)^{\mu} (\zj \xij)^q d(\zj\xij)^{\nu}-\sigma=0
\quad \left(\mbox{mod}\;
F^{i(p+\mu)+j(p+\mu+2q+2\nu)+1} \right).
\ea
\]

Statement~1 follows from
$\Big \langle\pi ,d(\zi\xij)d\left(\sum\limits_{k\ge j}z_k\xi_k
\right)\Big\rangle=0$, $i\ge j$.

Indeed, since
\[
\pi_{DS}=
\sum_i \zi\xii \dzi\wedge\dxi -
\sum_{i<j} \zi \zj  \dzi \wedge \dzj + \sum_{i<j} \xii \xij  \dxi \wedge
\dxj
+2\sum_{i>j} \zi \xii \dzj \wedge \dxj ,
\]
 we obtain
\[
\ba{l}
\ds \Big\langle \pi ,d(\zi\xij)d\left(\sum_{k\ge j}z_k\xi_k\right)\Big
\rangle=
\Big\langle \pi ,\sum_{k > i} \zi z_k d\xij d\xi_k
 + \sum_{k>i} \xij\xi_k dz_i dz_k
\vspace{3mm}\\
\ds \qquad -\zi\xij d\zj d\xij+\zi\xij d\zi d\xii +
 \sum_{k>i} (\zi\xi_k d\zj d\xi_k
 +\xij z_k d\zi d\xi_k)\Big\rangle
\vspace{3mm}\\
\ds \qquad =\sum_{k>j}\zi z_k\xij\xi_k +
\sum_{i>k\ge j}\zi z_k\xij\xi_k -
\sum_{k>i}\zi z_k\xij\xi_k
 - \zi\xij d\zj d\xij
\vspace{3mm}\\
\ds \qquad -2\sum_{k>j}\zi z_k\xij\xi_k +
 \zi\xij d\zi d\xii + 2\sum_{k>i}\zi z_k\xij\xi_k =0.
\ea
\]

Hence the Poisson homologies of Drinfeld--Sklyanin brackets are represented by
cocycles
\[
\sigma=(\zi \xij)^p d(\zi\xij)^{\mu} \left(\sum_{k\ge j}z_k \xi_k\right)^q
d\left(\sum_{k\ge j}z_k\xi_k\right)^{\nu}, \qquad i\ge j .
\]


Let us consider another example of computation for some standard
Poisson structures.

\medskip

\noindent
{\bf Example 3.}  Let
\[
 \pi=\sum_{i<j}\zi\zj\dzi\dzj.
\]

It is well-known that this Poisson structure arises from the skew-polynomial
deformation
\[
\zi\zj=q\zj\zi, \qquad i<j.
\]

It is quite easy to show that
\[
\dpi=\sum_{i<j}(\Lzi\izj-\izi\Lzj)=\sum_j
\left(\sum_{i<j}\Lzi-\sum_{j<i}\Lzi\right)\izj
\]
because $[\zi\dzi,\zj\dzj]=0$.

Now one can introduce the adjoint operator $\dpi^*$ and Laplacian $\Dpi$ as
follows:
\[
\dpi^*=\sum_j\left(\sum_{i<j}\Lzi-\sum_{j<i}\Lzi\right)\izj^*,
\]
where $\izi^*=d\zi\dzi$ and $\izi\izj^* +\izj^*\izi=\delta_{i,j}\Lzi$.

So
\[
 \Dpi =\dpi\dpi^*+\dpi^*\dpi
=\sum_j\left(\sum_{i<j}\Lzi-\sum_{j<i}\Lzi\right)^2 \Lzj.
\]
Recall that $\Lzi$ is grading operator, i.e. $\Lzi(\zj^k
d\zj^l)=(k+l)\delta_{i,j} \zi^k d\zi^l$.
Then, if $\omega$ is a ``harmonic'' homogeneous form and $m_i=deg_{\zi}\omega$,
we get the equation
\[
\sum_j\left(\sum_{i<j}m_i-\sum_{j<i}m_i\right)^2 m_j=0, \qquad m_j\ge 0.
\]
It means  that
$\left(\sum\limits_{i<j}m_i-\sum_{j<i}m_i\right) m_j=0$.

Let $ j_0=\max\{j,\; m_j\ne 0\}$. Then
$\left(\sum\limits_{i<j_0}m_i-\sum\limits_{j_0<i}m_i\right) m_{j_0}=0$
implies
$\sum\limits_{i<j_0}m_i=0$ or $m_i=0$, $i\ne j_0$.

Therefore, the solution is
\[
\bigcup_j \{(0,\ldots,m_j,0,\ldots, 0),\;m_j\in Z\}.
\]
Hence the space of ``$\Dpi$-harmonic'' forms, that naturally set
the space of homologies, is
\[
{\mathcal H}_{\pi} =\bigoplus_j C(\zj,d\zj)= \bigoplus_j C[\zj]\otimes C(d\zj).
\]


Now we consider the formal Poisson homology complex associated with the
Drinfeld--Sklyanin $r$-matrix structure.

\medskip

\noindent
{\bf Example 4}
In the af\/f\/ine coordinates on large Schubert cell of $CP^n$, the
Drinfeld--Sklyanin structure is as follows:
\[
\ba{l}
\ds \pi=\sum_{i<j} \zi\zj\dzi\wedge\dzj - \sum_{i<j} \bzi\bzj\bdzi\wedge\bdzj +
\sum_{i,j}(1+\delta_{i,j}) \zi\bzj\dzi\wedge\bdzj
\vspace{3mm}\\
\ds \qquad +2\sum_{i<j} \zi\bzi\dzj\wedge\bdzj
+2\left(\sum_k \zj\bzj\right)\sum_{i,j} \zi\bzj\dzi\wedge\bdzj.
\ea
\]
One can introduce a grading on the algebra of algebraic forms and polyvectors
with formal coef\/f\/icients such that all natural operations (i.e
exterior dif\/ferential and multiplication, interior multiplication and
Schouten--Nijenhuis brackets) preserve the grading. Here
we attach the grading degree $i$ to
$\zi$, $\bzi$, $d\zi$, $d\bzi$ and $-i$ to $\dzi$, $\bdzi$.
So the algebra of formal dif\/ferential forms $A=\Omega^*_{\mbox{\scriptsize form}}
(\zi,\bzi)$ is
graded. The corresponding increasing f\/iltration is
\[
F^p A=\bigoplus_{k\ge p}A_p, \qquad
A=F^0 A\supset F^1 A\supset \cdots \supset
F^p A\supset \cdots.
\]
The Poisson (Brylinsky) dif\/ferential $\dpi$ agrees with this
f\/iltration. In order to see it we decompose $\pi$ in two
components, $\pi=\pi_0+\pi_1$, where
\[
\ba{l}
\ds \pi_0=
\sum_{i<j} \zi\zj\dzi\wedge\dzj - \sum_{i<j} \bzi\bzj\bdzi\wedge\bdzj +
\sum_{i,j}(1+\delta_{i,j}) \zi\bzj\dzi\wedge\bdzj,
\vspace{3mm}\\
\ds \pi_1= 2\sum_{i<j} \zi\bzi\dzj\wedge\bdzj
+2\left(\sum_k \zj\bzj\right)\sum_{i,j} \zi\bzj\dzi\wedge\bdzj.
\ea
\]
The f\/irst component is of degree $0$, and the second one is the sum of bivector
f\/ields of positive weights. Thus,
$\dpi (F^p A)\subset F^p A$.
Moreover, in the corresponding spectral sequence $E^*_r$
\[
 E^*_0=\bigoplus_p F^pA / F^{p+1}A\simeq\bigoplus_p A_p= A , \qquad
 \delta_o= \delta_{\pi_0}. 
\]
 Thus, one has to compute the Poisson homology for $\pi_0$.
The Poisson dif\/ferential $\delta_{\pi_0}$ is
\[
\ba{l}
\ds \delta_{\pi_0}=\sum_{i<j}(\Lzi\izj-\izi\Lzj-\bLzi\bizj+\bizi\bLzj)
\vspace{3mm}\\
\ds \qquad +\sum_{i,j}(1+\delta_{i,j})(\Lzi\bizj-
\bizi\Lzj)
\vspace{3mm}\\
\ds \qquad =\sum_j\left(\sum_{i<j}\Lzi-\sum_{j<i}\Lzi +
\sum_i(1+\delta_{i,j})\bLzi\right)\izj
\vspace{3mm}\\
\ds \qquad
-\sum_j \left(\sum_{i<j}\bLzi-\sum_{j<i}\bLzi +\sum_i(1+\delta_{i,j})\Lzi
\right)\bizj.
\ea
\]
The homologies of complex are calculated as in the previous example:
\[
H_{\pi_0}(A) \simeq
{\cal H}_{\pi_0} =\bigoplus_j C(\zj,d\zj)\oplus \bigoplus_j C(\bzj,d\bzj).
\]
Since all ``harmonic'' forms are $\dpi$-closed, we see that the spectral
sequence converges in the f\/irst term and
\[
H_{\pi}(A) \simeq \bigoplus_j C(\zj,d\zj)\oplus \bigoplus_j C(\bzj,d\bzj).
\]

\subsection*{Acknowledgements}
The author is deeply indebted  to Vladimir Rubtsov  for
 helpful discussions.
I greatly acknowledge the warm  hospitality and support of the
Mathematics Department of the University of Angers where a part of
this results was obtained and the Theoretical Physics department
of Uppsala University. The last visit was possible due to the
grant INTAS 96-196. This work was partly supported by Russian
President's grant 96-15-96939 and RFBR grant 98-02-16575.

\label{kotov-lp}
\end{document}